\documentclass[12pt,a4paper]{article}
\usepackage{indentfirst}
\usepackage{latexsym}
\usepackage{graphicx}
\usepackage{mathrsfs}
\usepackage{amsmath,amsthm}
\usepackage{amssymb}
\usepackage{indentfirst}
\usepackage{mathtools}
\usepackage{cases}
\title{The proof of a conjecture about cages\thanks{This work was partially supported by the National Natural Science Foundation of China under Grant No. 19871034.}}
\author{Xiang-Feng Pan\thanks{Corresponding author: xfpan@ahu.edu.cn}, Jing-Zhong Mao, Hui-Qing Liu\\
Central China Normal University, Wuhan, 430079, P. R. China}

\date{Received August 8, 2000; Published April 25, 2001}

\begin{document}
\maketitle
\textbf{Remark: This paper was published in Chinese on April 25, 2001, in Volume 14, Issue 2 of the ``MATHEMATICA APPLICATA" on pages 99 to 102. This is the English version of the paper, incorporating additional details and rectifying certain typographical errors.}
\begin{abstract}
	The girth of a graph is defined as the length of a shortest cycle in the graph. A $(k; g)$-cage is a graph of minimum order among all $k$-regular graphs with girth $g$. A cycle $C$ in a graph $G$ is termed nonseparating if the graph $G-V(C)$ remains connected. A conjecture, proposed in [T. Jiang, D. Mubayi. Connectivity and Separating Sets of Cages. J. Graph Theory 29(1)(1998) 35--44], posits that every cycle of length $g$ within a $(k; g)$-cage is nonseparating. While the conjecture has been proven for even $g$ in the aforementioned work, this paper presents a proof demonstrating that the conjecture holds true for odd $g$ as well. Thus, the previously mentioned conjecture was proven to be true.
\end{abstract}
\section{Introduction}
\noindent Graphs considered in this paper are simple and finite. For undefined terminology and notations related to Graph Theory, refer to \cite{art1}.   

Let $G = (V(G), E(G))$ be a simple graph, where $V(G)$ is the set of vertices and $E(G)$ is the set of edges. The \textit{complement}  $\overline{G}$  of  $G$  is the simple graph whose vertex set is  $V(G)$  and whose edges are the pairs of nonadjacent vertices of  $G$. For vertices $u, v \in V(G)$, we denote the \textit{distance} between $u$ and $v$ in $G$ by $d_{G}(u, v)$; if $u$ and $v$ belong to different components of $G$, then we set $d_{G}(u, v) = \infty$.   The \textit{diameter} of $G$, denoted by $\operatorname{diam}(G)$, is defined as   
\[  
\max \{d_{G}(u, v) : \forall u, v \in V(G)\}.  
\]  
Let $C_n$ and $P_n$ denote a cycle and a path, respectively, each consisting of $n$ vertices. Given a path $P$ and a cycle $C$, their \textit{lengths}, denoted by $l(P)$ and $l(C)$, are defined as the number of edges in $P$ and $C$,  respectively.  A cycle of length $g$ is referred to as a \textit{$g$-cycle}.
If $G_{1}$ and $G_{2}$ are disjoint graphs, we denote their \textit{union} by $G_{1} + G_{2}$.   
For a subset $S \subseteq V(G)$, we define the \textit{neighborhood} of $S$ as   
\[  
N(S) = \bigcup_{v \in S} N(v).  
\]  
Furthermore, if $H$ is a subgraph of $G$, then the \textit{neighborhood of $S$ within $H$} is defined as   
\[  
N_{H}(S) = N(S) \cap V(H).  
\] 
The \textit{girth} of a graph $G$ is the length of a shortest cycle in $G$.\\
A $(k; g)$-\textit{graph} is a $ k $-regular graph with girth $g$. A $(k; g)$-\textit{cage} is a $(k; g)$-graph of minimum order.\\
A cycle $C$ in a graph $G$ is \textit{nonseparating} if $G-V(C)$ is connected.\\
The following conjecture is presented in \cite{art2}:\\

\noindent\textbf{Conjecture 1.} Every $g$-cycle in a $(k; g)$-cage is nonseparating.\\

\noindent\textbf{ Theorem 1\cite{art3}.} If $G$ is a simple graph with $|V(G)| \geq 3$ and $\delta \geq|V(G)| / 2$, then $G$ is Hamiltonian.\\

\noindent\textbf{Deﬁnition 1\cite{art2}.}
Suppose that $G$ is a $(k; g)$-cage and $H$ is an induced subgraph of $G$ with minimum degree $k-1$ in $H$. Let $B$ be the set of vertices of degree $k-1$ in $H$. For each permutation $\sigma$ of $B$, let 
\[  
D_{\sigma}(x, y) = d_{H}(x, y) + d_{H}(\sigma(x), \sigma(y)) + 2.  
\]  
We say that $H$ is a \textit{special} subgraph of $G$ if there exists a $\sigma$ such that $d_{H}(x, y) \geq \lceil g/4 \rceil-1$ and $D_{\sigma}(x, y) \geq g$ for all distinct $x,y\in B$.\\

\noindent\textbf{Lemma 1\cite{art2}.} If $ G $ is a $(k; g)$-cage and $ H $ is a special subgraph of $ G $, then $|V(H)| \geq|V(G)| / 2$.\\

\noindent\textbf{Lemma 2\cite{art2}.} If $ G $ is a graph with girth $g \geq 3$, $ S $ is a subset of $V(G)$ with $\operatorname{diam}(G[S])=d<g-2$, and $ H $ is subgraph of $G$ with $V(H) \cap S=\emptyset$, then every vertex of $N_{H}(S)$ has exactly one neighbor in $S$ and the distance in  $H$  between distinct neighbors of  $S$  is at least  $g-d-2$, with equality only if they are neighbors of a pair of vertices with distance  $d$  in  $G[S]$.\\

\noindent\textbf{Theorem 2\cite{art2}.}  If $k>3$ and $G$ is a $(k; g)$-cage, then $G$ is $3$-connected.\\

\noindent\textbf{Theorem 3\cite{art2}.} If $k \geq 3$ and $g \geq 4$ is even, then every $g$-cycle in a $(k; g)$-cage is nonseparating.\\

\noindent\textbf{Deﬁnition 2.} Let $G$ be a graph. We define a permutation $\sigma $ as a \textit{special permutation} of $V(G)$ if and only if for every edge $xy \in E(G)$, $\sigma(x)\sigma(y)$ is an edge in the complement $\overline{G}$ of $G$.\\

\noindent\textbf{Lemma 3.} If $C_g$ is a cycle of length $g \geq 5$, then there exists a special permutation of $V\left(C_{g}\right)$.\\
\begin{proof}  The following discussion will be divided into two cases.
	\begin{description}  
		\item[Case 1.] $g=5$. Label the vertices of $C_{5}$ clockwise as $1,2,3,4,5$. Then the cycle $135241$ is a Hamilton cycle of $\overline{C_{5}}$, denoted by $C$.  
		
		\item[Case 2.] $g>5$. It can be directly verified that $\overline{C_{g}}$ meets the criteria of Theorem 1, and thus there exists a Hamilton cycle in $\overline{C_{g}}$, which we also denote by $C$.
	\end{description}  
	
	For both cases, define a mapping $\sigma$ that maps the vertices in $C_{g}$ clockwise to the vertices in $C$ (also arranged clockwise). By definition, $\sigma$ is a special permutation on $V(C_{g})$.
\end{proof}  
\noindent\textbf{Lemma 4.} For any path $P_{n}$, $n \geq 4$, there exists a special permutation of $V\left(P_{n}\right)$.
\begin{proof}  
	Label the vertices of $P_{n}$ sequentially along $P_{n}$ as $1, 2, \ldots, n$.  We shall divide the proof into two distinct cases.
	
\begin{description}  
		\item[Case 1.] $n = 4$. Let $\sigma = (1342)$, where $(1342)$ is a cyclic permutation. It can be directly verified that $\sigma$ meets the requirements.  
		
		\item[Case 2.] $n > 4$. Denote by $P_{n}+ 1n$ a cycle of length $n$ whose vertices are labeled clockwise as $1,2,\ldots, n$. By Lemma 3, there must exist a special permutation $\sigma$ on $V(P_{n}+1n)$. By definition, $\sigma$ is also a special permutation on $V(P_{n})$.  
	\end{description}  
	\end{proof}  
\noindent\textbf{Lemma 5.} Let $G = C_{g} + P_{n}$ with $g > 5$ and $2\le n \le 3$. Then there exists a special permutation of $V(G)$.
\begin{proof} 
Denote the vertices of $C_{g}$ by $x_{1}, x_{2}, \ldots, x_{g}$ in a clockwise manner, and similarly, label the vertices of $P_{n}$ as $y_{1}, y_{2}, \ldots, y_{n}$ sequentially. By Lemma 4, there exists a special permutation $\sigma_{1}$ of $V(C_g)$.  

\begin{description}  
	\item[Case 1.]  $G=C_{g}+P_{2}$. Let\\
$$
\sigma(x)=\left\{\begin{array}{lll}
y_{1}, &x=x_{g}; \\
x_{g}, &x=y_{1}; \\
x, &x=y_{2} ; \\
\sigma_{1}(x), &x \in V(G)-\left\{x_{g}, y_{1}, y_{2}\right\} .
\end{array}\right.
$$
\item[Case 2.] $G=C_{g}+P_{3}$. Let\\
$$
\sigma(x)=\left\{\begin{array}{ll}
y_{2}, & x=x_{g}; \\
x_{g}, & x=y_{2}; \\
x, & x=y_{1} \text { or } y_{3} ; \\
\sigma_{1}(x), & x \in V(G)-\left\{x_{g}, y_{1}, y_{2}, y_{3}\right\}.
\end{array}\right.
$$.
\end{description}  
It is easy to show that $\sigma$ is a special permutation of $V(G)$ in both two cases.
\end{proof}
Let $G = \sum_{i=1}^{n} G_{i}$, $n > 1$, and let $\sigma_{i}$ be a permutation of $V(G_{i})$ for $i = 1, \ldots, n$. Let $\sigma_{1} \sigma_{2} \cdots \sigma_{n}$ denote the permutation $\sigma$ of $V(G)$ such that $\sigma(x) = \sigma_{i}(x)$ when $x \in V(G_{i})$.

\noindent\textbf{Lemma 6.} 
Let $P^{(i)}$ and $C^{(j)}$ denote paths and cycles, respectively, where $i = 1, 2, \ldots, t$ and $j = 1, 2, \ldots, q$. Suppose $l(C^{(j)}) \geq 5$ for each cycle $C^{(j)}$ and $l(P^{(i)}) \geq 1$ for each path $P^{(i)}$. Define the graph $G$ as  
\[  
G = \left( \sum_{i=1}^{t} P^{(i)} \right) + \left( \sum_{j=1}^{q} C^{(j)} \right).  
\]  
Furthermore, when $q = 0$ and $t = 1$, assume that the length of the single path $P^{(1)}$ is at least 3, i.e., $l(P^{(1)}) \geq 3$. Under these conditions, there exists a special permutation of the vertex set $V(G)$.

\begin{proof} If $q \neq 0$, then, by Lemma 3, for each $C^{(j)}$, where $j=1,2, \ldots, q$, there exists a special permutation $\tau_{j}$ of $ V(C^{(j)})$, respectively.\\
\begin{description}  
\item[Case 1.]  $\sum_{i=1}^{t}\left|V\left(P^{(i)}\right)\right| \geq 4$ and $ t=1 $. Then $\left|V\left(P^{(1)}\right)\right| \geq 4$. By Lemma
4,  there  must  exist  a  special  permutation $ \sigma_{1} $ of $V\left(P^{(1)}\right)$. Let $\sigma=\sigma_{1} \tau_{1} \tau_{2} \cdots \tau_{q}$, if $q>0 .$ Otherwise, let $\sigma=\sigma_{1}$.\\
\item[Case 2.]  $\sum_{i=1}^{t}|V(P^{(i)})| \geq 4$ and $t > 1$. Denote the start vertex and the end vertex of $P^{(i)}$ by $u_{i}$ and $v_{i}$, respectively. For all $i = 1, 2, \ldots, t-1$, let us connect $P^{(1)}, P^{(2)}, \ldots, P^{(t)}$ with an edge $v_{i}u_{i+1}$ between $P^{(i)}$ and $P^{(i+1)}$ to form a path $P$ of $\sum_{i=1}^{t}|V(P^{(i)})| \geq 4$ vertices. By Lemma 4, there exists a special permutation $\sigma_{1}$ of $V(P)$. Let $\sigma = \sigma_{1}\tau_{1}\tau_{2}\cdots\tau_{q}$ if $q > 0$; otherwise, let $\sigma = \sigma_{1}$.
\item[Case 3.]  $\sum_{i=1}^{t}\left|V\left(P^{(i)}\right)\right| \leq 3$ and $ t=1 $. Since  when $q=0$ and $t=1$ $l\left(P^{(1)}\right) \geq 3$, $q>0$. By Lemma 5, there exists a special permutation $ \sigma_{1} $ of $V\left(C^{(1)}+\left(\sum_{i=1}^{t} P^{(i)}\right)\right)$. Let $\sigma=\sigma_{1} \tau_{2} \tau_{3} \cdots \tau_{q,}$ if $q>2$. Otherwise, let $\sigma=\sigma_{1}$. 
\end{description}  
It is not diﬃcult to show that $\sigma$ is as required in all cases.
\end{proof}
Given a cycle $C$ of odd length $g$ and a vertex $u$ on $C$, there are two  vertices $v,w$ on $C$ with $d_{C}(u, v)=(g-1)/2$ and $d_{C}(u, w)=(g-1)/2$. We called $v,w$ the \textit{antipodal vertices} of $u$ on $C$, and $\{u,v\}$, $\{u,w\}$ the \textit{antipodal pairs} on $C$. 

\section{The proof of the conjecture}
\noindent \textbf {Theorem 4.} If $ k\ge 3$ and $g\ge 5$ is odd, then every $ g $-cycle in a $(k; g)$-cage is nonseparating.
\begin{proof} 
Among the set $AC$ of $g$-cycles in $G$ such that $G-V(C)$ is disconnected, select $C\in AC$ to minimize the order of the smallest component $H$ of $G - V(C)$.\\
Since $\operatorname{diam}(C)=(g-1) / 2<g-2$ (when $g\ge 5$), by Lemma 2, every vertex $ w $ in $ N_{H}(C)$ has a unique neighbor, which we denote by $w^{\prime}$, in $C$.  Hence, $H$ has minimum degree $ k-1 $, with $ N_{H}(C) $ being the set of vertices with degree $k-1$. 
For all $x, y \in N_{H}(C)$, a shortest $(x,y)$-path in $H$, along with a shortest $(x^{\prime}, y^{\prime})$-path in $C$, and the edges $xx^{\prime}, yy^{\prime}$, collectively form a cycle in $G$ with a length of at least $g$.
Thus
$$
d_{H}(x, y) \geq g-2-d_{C}\left(x^{\prime}, y^{\prime}\right).
$$
Since $d_{C}\left(x^{\prime}, y^{\prime}\right) \leq(g-1) / 2$,  we have $d_{H}(x, y) \geq(g+1) / 2-2$, with equality only if $d_{C}\left(x^{\prime}, y^{\prime}\right)=(g-1) / 2$. We shall call a pair $\{x, y\}$ a \textit{bad pair} in $N_{H}(C)$, if $d_{H}(x, y)=(g+1) / 2-2$. As we have seen, if $\{x, y\}$ is a bad pair, then $\left\{x^{\prime}, y^{\prime}\right\}$ is an antipodal pair on $ C $. A pair $\{x, y\}$ in $N_{H}(C)$ that is not a bad pair satisﬁes $d_{H}(x, y) \geq(g+1) / 2-1$.\\
If we can demonstrate the existence of a permutation $\sigma$ of $N_{H}(C)$ such that $\{\sigma(x), \sigma(y)\}$ is not a bad pair for any bad pair $\{x, y\}$ within $N_{H}(C)$, then $\sigma$ satisfies the following condition:   
$D_{\sigma}(x, y) = d_{H}(x, y) + d_{H}(\sigma(x), \sigma(y)) + 2 \geq (\frac{g+1}{2} - 2) + (\frac{g+1}{2} - 1) + 2 = g$,  
for every pair $x, y \in N_{H}(C)$. Consequently, by invoking Lemma 1, we can infer that $|V(H)| \geq \frac{|V(G)|}{2}$, thereby obtaining a contradiction.\\
Let us ﬁrst study all bad pairs in $N_{H}(C)$. Let $A=\{x|\{x, y\}$ is a bad pair in $N_{H}(C)\}$, $B=\left\{xy|\{x, y\}\right.$ is a bad pair in $\left.N_{H}(C)\right\}$. Let $G_{1}$ be the graph with $V\left(G_{1}\right)=A$ and $E\left(G_{1}\right)=B$.\\
\textbf{Claim 1.} If $\{x, y\},\{x, z\}$ are two bad pairs, then $y^{\prime} \neq z^{\prime}$.\\
If $y^{\prime}=z^{\prime}$, then a shortest $(x, y)$-path and a shortest $(x, z)$-path in $H$, together with the edges $zz^{\prime}$ and $yy^{\prime}$, form a cycle in $G$ of length at most $g-1$, which contradicts the fact that $G$ has girth $g$.\hfill$\qed$\\
\textbf{Remarks.} Since there are at most two antipodal vertices of $x^{\prime}$ on $C$, by Claim 1, it follows that $1\le d_{G_{1}}(x) \leq 2$ for all $x \in N_{H}(C)$. 
Furthermore, we are cognizant that $G_{1}$ constitutes a disjoint union of cycles and paths, wherein each path, if present, possesses length at least $1$.\\
\textbf{Claim 2.} If there is a cycle component of $G_{1}$, then all cycles in $ G_{1} $ are of length at least $g$.\\
Let $C^{*}$ be a cycle component of $G_{1}$ with length $m$. Denote the vertices of $C^{*}$ by $x_{1}, x_{2}, \ldots, x_{m}$ clockwise one by one. Let $C^{+} = C + \sum_{i=1}^{m} x_{i}^{\prime} x_{i+1}^{\prime}$ (admitting repeated edges), where $i+1$ is taken modulo $m$. Then $x^{\prime}y^{\prime} \in E\left(C^{+}\right)$ for every antipodal pair $\left\{x^{\prime}, y^{\prime}\right\}$ of $C$. Otherwise, assume $\left\{z^{\prime}, u^{\prime}\right\}, \left\{z^{\prime}, w^{\prime}\right\}$ are two distinct antipodal pairs of $C$ such that $z^{\prime}u^{\prime}\in E\left(C^{+}\right)$ but $z^{\prime}w^{\prime} \notin E\left(C^{+}\right)$. Let $z^{\prime}$ be a neighbor of $z$ in $G$ on $C$. Then, by Claim 1, $d_{G_{1}}(z) = 1$, which contradicts the fact that $C^{*}$ is a cycle component of $G_{1}$.
 Thus, $m = \left|E\left(C^{+}\right)\right| - |E(C)| \geqslant$ the total count of antipodal pairs in $C$, which is equal to $g$.\hfill$\qed$\\
\textbf{Claim 3.} If $G_{1}$ is a path $P$, then $l(P) > 3$; otherwise, $\exists w \in N_{H}(C) - A$.\\
In fact, assume that $l(P)<3$ and $N_{H}(C)-A=\emptyset$. By the definition of $A, A \subseteq N_{H}(C)$. Thus $N_{H}(C)=A$. There are two cases:\\
$i$) $l(P)=1$. Let $P=uv$. Since $G_{1}$ is only a path $P, N_{H}(C)=A=\{u, v\}$. Thus $\left|N_{H}(C)\right|=2$. Because $G-N_{H}(C)$ is not connected, $ N_{H}(C) $ is a cut  set of $ G $. However, by Theorem 2, $ G $ is 3-connected, a contradiction.\\
$ii$) $l(P)=2$. Let $P=u v w$. Then $N_{H}(C)=A=\{u, v, w\}$. By the deﬁnition of $G_{1},\{u, v\},\{v, w\}$ are two bad pairs in $N_{H}(C)$. Therefore, $d_{H}(u, v)=d_{H}(v, w)=(g+1) / 2-2$ and $\left\{u^{\prime}, v^{\prime}\right\},\left\{v^{\prime}, w^{\prime}\right\}$ are antipodal pairs on $C$. By Claim 1, $u^{\prime} \neq w^{\prime}$. So by the deﬁnition of antipodal pair, $u^{\prime} w^{\prime} \in E(C)$. Observe a cycle $ C^{0} $ of length $ g $ which consists of a shortest $ uv $-path and a shortest $ vw $-path in $ H $, and three edges $u u^{\prime}, w w^{\prime}, u^{\prime} w^{\prime}$. Since $N_{H}(C)=A=\{u, v, w\}$ and $k \geq 3$, $G-V(C^{0})$ is not connected and there exists a component $ H_{0} $ of $G-V(C^{0})$ with $H_{0} \subseteq H-\{u, v, w\}$ but $\left|V\left(H_{0}\right)\right| \leq|V(H)-\{u, v, w\}|<|V(H)|$, contrary to the choice of $ H $.\hfill$\qed$\\
Now we try to ﬁnd a permutation $ \sigma $ of $ N_{H}(C) $ such that $\{\sigma(x), \sigma(y)\}$ is not a bad pair for every bad pair $\{x, y\}$ in $N_{H}(C)$.\\
\textbf{Case 1.} $G_{1}$ is only a path $P$ with $l(P)<3$. By Claim 3, there exists a vertex $w \in N_{H}(C)-A$.\\
\textbf{Subcase 1.1.} $l(P)=1$. Let $\{u, v\}$ be the unique bad pair in $N_{H}(C)$. Let
$$
\sigma(x)=\left\{\begin{array}{ll}
w, & x=u; \\
u, & x=w;\\
x, & x \in N_{H}\left(C\right)-u-w.
\end{array}\right.
$$
\textbf{Subcase 1.2.} $l(P)=2$. Let $\{u, v\},\{v, z\}$ be all bad pairs in $N_{H}(C)$. Let
$$
\sigma(x)=\left\{\begin{array}{ll}
w, & x=v; \\
v, & x=w; \\
x, & x \in N_{H}(C)-v-w.
\end{array}\right.
$$
\textbf{Case 2.} $G_{1}=\left(\sum_{i=1}^{l} P^{(i)}\right)+\left(\sum_{j=1}^{q} C^{(j)}\right)$, where $P^{(i)}$ denote paths and $C^{(j)}$ cycles of length at least $ g $, and when $q=0$ and $l=1, l\left(P^{(1)}\right) \geq 3 .$ By Lemma 6, there exists a special permutation $ \sigma_{1} $ of $ V(G_{1}) $. Let 
$$
\sigma(x)=\left\{\begin{array}{ll}
\sigma_{1}(x), & x \in A; \\
x, & x \in N_{H}(C)-A.
\end{array}\right.
$$
It is straightforward to demonstrate that $\sigma$ is the permutation we aim to obtain in all cases, thereby concluding the proof of Theorem 4.
\end{proof}
\section{Conclusion} By Theorems 3 and 4, Conjecture 1 is true.


\begin{thebibliography}{99}
	\bibitem{art1}J.A. Bondy, U.S.R. Murty. Graph Theory with Application. The Macmillan Press LTD, 1976.
	\bibitem{art2}T. Jiang, D. Mubayi. Connectivity and Separating Sets of Cages. J. Graph Theory 29(1998) 35-44.
	\bibitem{art3}G.A. Dirac. Some Theorems on Abstract Graphs. Proc. London Math. Soc. 2 (1952) 69-81.
\end{thebibliography}
\end{document}